

%
%
%
%
%
%
\magnification\magstep1
\baselineskip = 18pt
\def\n{\noindent}
\overfullrule = 0pt
\def\n{\noindent}  
 
\def\qed{{\hfill{\vrule height7pt width7pt
depth0pt}\par\bigskip}}

\def\pf{ {\n {\bf Proof.}\ }}

\centerline{\bf A simple proof of a theorem}
\centerline{\bf of Kirchberg and related results on $C^*$-norms}\bigskip
\centerline{by}
\centerline{Gilles Pisier\footnote*{Partially supported by
the NSF}}\bigskip

\n {\bf \S 0. Introduction. Main results.}

Recently, E.\ Kirchberg [K1--2] revived the study of pairs of
$C^*$-algebras $A,B$ such that there is only one $C^*$-norm on the
algebraic tensor product $A\otimes B$, or equivalently such that $A
\otimes_{\rm min}B = A\otimes_{\rm max}B$. Recall that a $C^*$-algebra is
called nuclear cf.\ [L, EL] if this happens for any $C^*$-algebra $B$.
Kirchberg [K1] constructed the first example of a non-nuclear $C^*$-algebra
such that $A\otimes_{\rm min} A^{op} = A \otimes_{\rm max} A^{op}$. He also
proved the following striking result [K2] for which we
give a very simple proof and which we extend.

\proclaim Theorem 0.1. {\bf (Kirchberg [K2]).} Let $F$ be
any free group and let $C^*(F)$ be the (full)
$C^*$-algebra of $F$, then $$C^*(F) \otimes_{\rm min} B(H)
= C^*(F) \otimes_{\rm max} B(H).$$

Here, and throughout the paper 
(unless specified otherwise) $H$ is a separable infinite
dimensional Hilbert space
and $B(H)$ is the space of all bounded operators on $H$.

More generally, we will prove

\proclaim Theorem 0.2. Let $(A_i)_{i\in I}$ be a family of unital
$C^*$-algebras such that
$$A_i\otimes_{\rm min} B(H) = A_i\otimes_{\rm max} B(H)\leqno \forall i\in
I$$
then the free product $A = \mathop{*}\limits_{i\in I}A_i$ (in the category
of unital $C^*$-algebras) satisfies
$$A \otimes_{\rm min} B(H) = A\otimes_{\rm max} B(H).$$

\proclaim Corollary 0.3. Let $(G_i)_{i\in I}$ be a family
of discrete amenable groups and let $G=\mathop{*}\limits_{i\in
I}G_i$ be their free product. Then $$C^*(G) \otimes_{\rm
min} B(H) = C^*(G) \otimes_{\rm max} B(H).$$

Our method strongly uses the theory of operator spaces, recently developed
in a series of papers by Effros-Ruan [ER1-2] and
Blecher-Paulsen [BP]. A key observation is that when a
$C^*$-algebra is generated by a finite set of unitaries,
then the operator space structure of the linear span of
these generators (up to complete isometry) is enough to
encode some important properties of the $C^*$-algebra they
generate.

 {\bf Notation:}\ Our notation is standard, except 
perhaps that we denote by
$E_1\otimes E_2$ the algebraic tensor product of two vector spaces.

 {\bf Acknowledgement:}\ I am grateful 
to Uffe Haagerup and to Eberhard
Kirchberg for useful conversations concerning
Proposition~6, during two different visits at the Fields
Institute in Waterloo in the spring of 1995.
Thanks also to F. Boca
for a reference, to S. Wassermann for a
correction and to S. Popa for a useful comment.

\n {\bf \S 1. Proofs.}

The main idea of our proof of Theorem~0.1 is that if $E$ is the linear span
of 1 and the free unitary generators of $C^*(F)$, then it suffices to check
that the min- and max-norms coincide on $E\otimes B(H)$. More generally, we
will prove

\proclaim Theorem 1. Let $A_1, A_2$ be unital $C^*$-algebras. Let
$(u_i)_{i\in I}$ (resp.\ $(v_j)_{j\in J}$) be a family of unitary operators
which generate $A_1$ (resp.\ $A_2$). Let $E_1$ (resp.\ $E_2$) be the closed
span of $(u_i)_{i\in I}$ (resp.\ $(v_j)_{j\in J}$). Assume  $1\in E_1$ and
$1\in E_2$. Then the following assertions are equivalent:\medskip
\item{(i)} The inclusion map $E_1\otimes_{\rm min}E_2\to A_1\otimes_{\rm
max} A_2$ is completely isometric.
\item{(ii)} $A_1\otimes_{\rm min}A_2 = A_1\otimes_{\rm max} A_2$.
\medskip

{ We  will use several elementary facts, as follows.
The first one is a well known property of unitary
dilations.}

\proclaim Lemma 2. Let $u\in B({\cal H}), \hat u\in B(\widehat H)$ be
unitaries and let $S\colon \ {\cal H}\to \widehat H$ be an isometry with
range $K\subset \widehat H$, such that
$$u = S^*\hat uS.$$
Then $K  = S({\cal H})$ is invariant under $\hat u$ and $\hat u^*$, so that
$\hat u$ commutes with $P_K$.

\n {\bf Proof.} Note that $P_K = SS^*$. We have $\forall h
\in {\cal H}\ \|P_K\hat uSh\| = \|S^*\hat uSh\|$ hence
$$\|\hat u S(h)\|^2 = \|S^*\hat u Sh\|^2 + \|(1-P_K) \hat uSh\|^2$$
hence
$$\|(1-P_K)\hat u Sh\|^2 = \|\hat uSh\|^2 - \|S^*\hat
uSh\|^2$$ so that if $S^*\hat uS$ is isometric this is
$= 0$. Taking adjoints we obtain the same for  $\hat
u^*$. \qed

\proclaim Lemma 3. Let $(a_i)_{i\in I}$ and $(b_i)_{i\in
I}$ be finitely supported families of operators in
$B({\cal H})$. We have $$\left\|\sum_{i\in
I}a_ib_i\right\| \le \left\|\sum a_ia^*_i\right\|^{1/2}
\left\|\sum b^*_ib_i\right\|^{1/2}.\leqno (1)$$

\n {\bf Proof.} This is an easy consequence of the Cauchy-Schwarz
inequality and is left to the reader.\qed

\n {\bf Remark.} For any unitary $u_i$, the norm $\big\|\sum
a_iu_ib_i\big\|$ is $\le$ the right side of (1). However, the norm of
$\sum\limits_{i\in I} a_ib_iu_i$ can be larger in general.

The next result is known, I personally
 learned this kind of
result from Paulsen. (see e.g. [Pa3].)

\proclaim Lemma 4. Let $F$ be a free group. Let $(U_i)_{i\in I}$ be the
free unitary generators of $C^*(F)$ (i.e.\ these are the
unitaries
corresponding to the free generators of
$F$ in the universal 
unitary representation of $F$). Let $(x_i)_{i\in
I}$ be a finitely supported family in $B(H)$. Then the
following are equivalent\medskip \item{\rm (i)} The linear
map $T\colon \ \ell_\infty(I) \to B(H)$ defined by
$T((\alpha_i)_{i\in I}) = \sum\limits_{i\in I}
\alpha_ix_i$, satisfies $$\|T\|_{cb}<1.$$
 \item{\rm (ii)} We have $$\left\|\sum_{i\in I}U_i\otimes
x_i\right\|_{C^*(F)\otimes_{\rm min} B(H)} <1.$$
\item{(iii)} There are operators $a_i,b_i$ in $B(H)$ such
that $x_i=a_ib_i$ and
$$\left\|\sum a_ia^*_i\right\|^{1/2} \left\|\sum b^*_ib_i\right\|^{1/2} <
1.$$
Moreover, the same remains true if we add the unit element
to the family $(U_i)_{i\in I}$.

\n {\bf Proof.} It is easy to check 
going back to the definitions
of both sides
that $$\|T\|_{cb} =
\big\|\sum U_i\otimes x_i\big\|_{\rm min}.$$
We leave this as an exercise for the reader.
This shows the equivalence of (i) and (ii).

\n Now assume (i). 
 By the
factorization of $c.b.$ maps we can write $T(\alpha) =
V^*\pi(\alpha)W$ where $\pi\colon \ \ell_\infty(I) \to
B(\widehat H)$ is a representation and where $V,W$ are in
$B(H, \widehat H)$ with $\|V\|\, \|W\| = \|T\|_{cb}$. 
 We can assume $I$ finite and $\widehat H=H$. Let
$(e_i)_{i\in I}$ be the canonical basis of $ \ell_\infty(I)$,
 we set $$a_i =
V^*\pi(e_i)\quad \hbox{and}\quad b_i = \pi(e_i)W.$$ It is
then easy to check (iii). Finally, the implication
(iii)~$\Rightarrow$~(ii) follows from Lemma~3
(applied to $U_i\otimes a_i$ and $1\otimes b_i$).
The last assertion follows from the next remark.\qed

\n {\bf Remark 5.} Let $A$ be a $C^*$-algebra and let
$(a_i)_{i\in I}$ be a finitely supported family of
elements of $A$. Let $F$ be a  free group freely generated
by a family $(g_i)_{i\in I}$ and let $U_i$ be the
associated unitaries in $C^*(F)$. Fix $i_0\in I$ and let
$I'=I - \{i_0\}$. We wish to verify here that for if 
$\alpha$ is either the minimal
or the maximal $C^*$-norm on $C^*(F)\otimes A$ we have
$$\alpha\left(1 \otimes a_{i_0} + \sum_{i\in I'}
U_i\otimes a_i\right) = \alpha \left(\sum_{i\in I} U_i
\otimes a_i\right).\leqno (2)$$ This can be justified as
follows.(This kind of result was also pointed out to me by
Vern Paulsen.)

 Consider the
family $(\gamma_i)_{i\in I}$ in $F$ defined as follows
$$\gamma_i = g^{-1}_{i_0}g_i\qquad \forall i\in I' \quad \hbox{and}\quad
\gamma_{i_0} = g_{i_0}.$$
We claim that $(\gamma_i)_{i\in I}$ are again free in $F$ and generate $F$.
This is easy and left to the reader. It follows that the map which takes
each $g_i$ to $\gamma_i$ extends to an automorphism $h\colon \ F\to F$.
This automorphism $h$ induces an isometric unital
representation $\pi\colon\ C^*(F)\to C^*(F)$ which takes
$U_i$ to $U^*_{i_0}U_i$ for all $i$ in $I'$. Now let
$L\colon \ C^*(F)\to C^*(F)$ be the operation of left
multiplication by $U_{i_0}$. Then the composition $L\pi
\otimes I_A$ clearly is isometric with respect to the
minimal
or the maximal $C^*$-norm but it preserves $U_i \otimes
a_i$ for all $i$ in $I'$ and takes $1\otimes a_{i_0}$ to
$U_{i_0} \otimes a_{i_0}$. This yields (2).\qed

The following simple fact is essential in our argument.

\proclaim Proposition 6. 
Let $A,B$ be two unital $C^*$-algebras.
Let $(u_i)_{i\in I}$ be a family of unitary elements
of $A$ generating $A$ as a unital $C^*$-algebra
(i.e.\ the smallest unital $C^*$-subalgebra of $A$
containing them is $A$ itself).
Let $E\subset A$ be
the  linear span of $(u_i)_{i\in I}$
and $1_A$.
 Let $T\colon\ E\to B$
be a linear operator
such that $T(1_A)=1_B$ and taking each
$u_i$ to a unitary in $B$.
Then, $\|T\|_{cb}
\le 1$ suffices to ensure that $T$ extends to a
(completely) contractive representation
(=$*$-homomorphism) from $A$ to  $B$.

\n {\bf Proof.} 
Consider $B$ as embedded into $B({\cal H})$. Then,
it clearly suffices to prove this statement for 
$B=B({\cal H})$, which we now assume.
By the Arveson-Wittstock extension theorem
(cf.\  [Pa1, p. 100]), $T$ extends to
a
complete contraction $\hat T\colon \ A\to B({\cal H})$.
Since $T$ is assumed
unital, $\hat T$ is unital, hence by a well known result
(cf.\ e.g.\ [Pa1]) $\hat T$ must be completely positive,
and more precisely (cf.\ [Pa1]) we can write
$$\hat T(x) = S^*\hat \pi(x)S$$
where $\hat\pi\colon \ A\to B(\widehat H)$ is 
a unital representation (=$*$-homomorphism) and
$S\colon \ {\cal H}\to \widehat H$ is an isometry. Now for
 any unitary $U$ in the family $(u_i)_{i\in I}$,
 we have
$T(U) = \hat T(U) = S^*\hat \pi(U)S$, hence by Lemma~2,
since $T(U)$ is unitary by assumption, if
$K=S({\cal H})$ then $P_K$ commutes with $\hat \pi(U)$.
Now since these operators $U$ generate $A$, this implies
that $P_K$ commutes with $\hat \pi(A)$, so that $\hat T$
is actually a $*$-homomorphism. Thus, $\hat T$ is an
extension of $T$ and a (contractive) $*$-homomorphism.
This completes the proof. [Note that, a posteriori,
the Arveson-Wittstock 
completely contractive extension $\hat T$ is unique.  In
short, the proof reduces to this:\ the multiplicative
domain of $\hat T$  is a unital $C^*$-algebra (cf.\
[Ch1]) and contains $(u_i)_{i\in I}$, hence it
is equal to $A$.]\qed

\n {\bf Remark.} We will apply Proposition 6 in the
following particular situation.
 Let ${\cal A}\subset A$ be the (dense)
unital $*$-algebra generated by $E$. Consider a unital
$*$-homomorphism $u\colon \ {\cal A}\to B$. Then
$\|u_{|E}\|_{cb} \le 1$ suffices to ensure that $u$
extends to a (completely) contractive representation
(=$*$-homomorphism) from the whole of $A$ to  $B$.

\n {\bf Remark.} In the same situation as in Proposition
6, note that, if $T$ is a complete isometry,
then  $\hat T$ is a faithful representation
onto the $C^*$-algebra $B_1$ generated by the range of
$T$. Indeed, by Proposition 6 applied to $T^{-1}$,  $\hat
T\colon \ A\to B_1$ is left invertible. This can be used
to give
a very simple proof of the fact due
to Choi ([Ch2]) that the full $C^*$-algebra of any free
group admits a faithful representation into a direct sum
of matrix algebras. By Proposition 6, it suffices
to check this on the free generators and this is quite
easy.

\n {\bf Proof of Theorem 1.} The implication (ii) $\Rightarrow$ (i) is
trivial, so we prove only the converse. Assume (i). 
Let $E= E_1\otimes_{\rm min}  E_2$. We view $E$ as a subspace
of $A= A_1 \otimes_{\rm min} A_2$.
By (i), we have an  inclusion map
  $T\colon\ E_1\otimes_{\rm min}  E_2\to A_1
\otimes_{\rm max} A_2$ with $\|T\|_{cb}\le 1$.
By Proposition 6, $T$ extends to
a (contractive) representation $\hat T$ 
from $A_1 \otimes_{\rm min} A_2$
to $A_1 \otimes_{\rm max} A_2$.
Clearly $\hat T$ must preserve 
 the algebraic
tensor products $A_1 \otimes 1$
and $1 \otimes A_2$, hence also
 $A_1 \otimes A_2$. Thus we obtain (ii).\qed

\n {\bf Remark.} Let us denote
by $E_1\otimes 1 +1\otimes E_2$ the  linear subspace
spanned by elements of $A_1 \otimes A_2$ of the form
$\{a_1\otimes 1 + 1 \otimes a_2\}$. Then, in the situation
of Theorem 1, $E_1\otimes 1+1\otimes E_2$ generates
$A_1\otimes_{\rm min} A_2$, so that it suffices for the conclusion
of Theorem 1 to assume that the
operator space structures induced
on $E_1\otimes 1+1\otimes E_2$
by the min and max norms
  coincide.

\n {\bf Proof of  Kirchberg's Theorem 0.1.} Let $A_1 =
C^*(F),\  A_2=B(H)$. We take $E_2=B(H)$ and let $E_1$ be
the linear span of the unit and the free unitary generators $(U_i\mid
i\in I)$ of $C^*(F)$ (i.e.\  associated to  the free
generators of $F$).

\n Consider $x\in E_1\otimes E_2$, with $\|x\|_{\rm
min}<1$. By Lemma~4 we can write $x = \sum\limits_{i\in I}
U_i\otimes x_i$ with $x_i \in B(H)$, $(x_i)_{i\in I}$
finitely supported, admitting a decomposition as $x_i =
a_ib_i$ with $\big\|\sum a_ia^*_i\big\|<1$, $\big\|\sum
b^*_ib_i\big\| <1$, $a_i,b_i\in B(H)$. Now, let $\pi\colon
\ A_1\otimes_{\rm max} A_2\to B({\cal H})$ be any faithful
representation. Let $\pi_1 = \pi_{|A_1\otimes 1}$ and
$\pi_2 = \pi_{|1\otimes A_2}$. We have $$\eqalign{\pi(x)
&=\sum_{i\in I} \pi_1(U_i) \pi_2(x_i)\cr &= \sum_{i\in I}
\pi_1(U_i) \pi_2(a_i) \pi_2(b_i)}$$ hence, since $\pi_1$
and $\pi_2$ have commuting ranges we have $\pi(x) = y$
with $$y = \sum_{i\in I} \pi_2(a_i) \pi_1(U_i)
\pi_2(b_i).$$ Now by Lemma 3 (and the remark following it)
we have $$\|y\|\le\left\|\sum_{i\in I} \pi_2(a_i)
\pi_2(a_i)^*\right\|^{1/2} \left\|\sum \pi_2(b_i)^*
\pi_2(b_i)\right\|^{1/2} < 1.$$ Hence we conclude that
$$\|x\|_{\rm max} = \|\pi(x)\| < 1.$$
This shows that the min and max norms coincide
on $E_1\otimes B(H)$, but since
$M_n(B(H))\approx B(H)$ for any $n$, this implies
``automatically" that the inclusion
$$E_1\otimes_{\min} B(H)\to A_1\otimes_{\max} B(H)$$
is completely isometric. In other words,
the operator space structures
associated to the min and max norms coincide.
Thus, we conclude by Theorem 1. 
{\qed}

We can prove the following extension of Kirchberg's theorem.

\proclaim Theorem 7. Let $(A_i)_{i\in I}$ be a 
family of 
$C^*$-algebras (resp. unital $C^*$-algebras). Assume that
for each $i$ in $I$ $$A_i \otimes_{\rm min} B(H) = A_i
\otimes_{\rm max} B(H).\leqno (3)$$ 
We will denote by ${\mathop{\dot{*}}\limits_{i\in I} A_i}$
(resp. by $\mathop{*}\limits_{i\in I} A_i$)
  their free product
in the
category of  $C^*$-algebras (resp.   in the
category of unital $C^*$-algebras). Then we have
$$\left({\mathop{\dot{*}}\limits_{i\in I} A_i}\right)
\otimes_{\rm  min} B(H) =
 \left({\mathop{\dot{*}}\limits_{i\in I} A_i}\right)
\otimes_{\rm max} B(H),\leqno {(\dot{4})}$$ and in the
unital case $$\left(\mathop{*}\limits_{i\in I} A_i\right)
\otimes_{\rm  min} B(H) = \left(\mathop{*}\limits_{i\in I}
A_i\right) \otimes_{\rm max} B(H).\leqno (4)$$

\n {\bf Remark.} Kirchberg's theorem for $F=F_I$ corresponds to $A_i =
C^*({\bf Z})$ for all $i$ in $I$ (in the unital case).

The next result is well known, by now. It
is a corollary of the Paulsen-Smith extension of the Christensen-Sinclair
factorization of bilinear maps. 
We only sketch the standard argument.
We denote by $(y_1,y_2) \longrightarrow y_1 \odot y_2$
the natural bilinear map from $A_1 \otimes B(H) \times A_2
\otimes B(H)$ to $(A_1\otimes A_2) \otimes B(H)$ which
takes $$(a_1\otimes b_1, a_2\otimes b_2)\quad
\hbox{to}\quad a_1 \otimes a_2 \otimes b_1b_2.$$

\proclaim Lemma 8. Let $A_1,A_2$ be two operator spaces
and let $y \in A_1 \otimes A_2 \otimes B(H)$, with $H$
infinite dimensional. Then $\|y\|_{(A_1\otimes _h
A_2)\otimes_{\rm min}B(H)} <1$ iff there is a factorization
$$y=y_1 \odot y_2$$
with $y_i \in A_i\otimes B(H)$ such that $\|y_i\|_{\rm min}<1$ for $i=1,2$.

\n {\bf Proof.} (Sketch)\ We may assume $A_1,A_2$ both finite dimensional.
Then, by the self duality (and other properties) of the
Haagerup tensor product (cf. \ [ER2] and [BP]), $y$
can be identified with a linear map $\tilde y\colon \
A^*_1\otimes_h A^*_2\to B(H)$ with $\|\tilde y\|_{cb}<1$.

\n By the
factorization theorem (cf.\ [PS]) we have maps $\sigma_i\colon \ A^*_i \to
B(\widehat H)$ with $\|\sigma_i\|_{cb}<1$ and operators $V,W$ in the unit
ball of $B(H,\widehat H)$ such that $\tilde y(\xi_1\otimes \xi_2) =
W^*\sigma_1(\xi_1)\sigma_2(\xi_2)V$. But since $H$ is infinite dimensional,
we can assume $\widehat H=H$ and (incorporating $W^*$ and $V$ into
$\sigma_1$ and $\sigma_2)$ we can get rid of $W^*$ and $V$. Thus we obtain
$\tilde y_i\colon \ A^*_i\to B(H)$ with $\|\tilde y_i\|_{cb}<1$ such that
$$\tilde y(\xi_1\otimes \xi_2) = \tilde y_1(\xi_1) \tilde
y_2(\xi_2)\quad \forall \xi_i\in A_i^*. \leqno (5)$$
Returning to the tensor products, $\tilde y_i$
corresponds to $y_i\in A_i \otimes B(H)$ with
$\|y_i\|_{\rm min}<1$ and (5) means that $$y=y_1\odot
y_2.$$ 
This yields the desired factorization.
{\qed}

The key point is the following important observation concerning the
Haagerup tensor product (part (ii) in Lemma~9 is perhaps the main new idea
of this paper).

\proclaim Lemma 9. Let $A_1,A_2$ be two  $C^*$-algebras
(resp. unital)  satisfying (3).
Let $A_1\dot{*}A_2$ (resp. $A_1{*}A_2$) be their free
product (resp. free product as unital $C^*$-algebras) and
let $E \subset A_1\dot{*}A_2$ (resp. $E \subset A_1*A_2$)
be the linear span in $A_1\dot{*}A_2$ (resp. $A_1*A_2$) of
all elements of the form $a_1a_2$ with $a_i\in A_i$. Then:
\medskip \item{\rm (i)} The mapping $p\colon \ a_1\otimes
a_2 \to a_1a_2$ is a complete  isometry of $A_1\otimes_h
A_2$ onto the closure of $E$ in $A_1\dot{*}A_2$ (resp. $A_1{*}A_2$). 
\item{\rm (ii)} The min
and max norms of $(A_1\dot{*}A_2) \otimes B(H)$ (resp.
$(A_1*A_2) \otimes B(H)$) coincide on $E\otimes
B(H)$.\medskip

\n {\bf Proof.} The assertion (i) is essentially known. It is proved in
[CES] for the non-unital free product, and nothing is
said there about the unital case.
However, when
 $A_1,A_2$ are unital, the argument of [CES]
can  be pursued to yield (i) as stated
above. A similar argument is used in [Ha]. As far as we
know,  this question is nowhere considered
(except for [Ha]). Therefore, we decided to include the
details: \n we will now verify (i) in the unital case,
starting from the non-unital case, which is treated in
[CES].

\n Consider $x = \sum a^1_i \otimes a^2_i$
 in $A_1\otimes A_2$. By [CES] we
have
$$\|x\|_h = \sup\left\{\left\| \sum \sigma_1(a^1_i) \sigma_2(a^2_i)\right\|
\right\}\leqno(6)$$
where the supremum runs over all pairs $\sigma_i\colon \ A_i \to B(\widehat
H)$ of (not necessarily unital) $*$-homomorphisms, with
$\widehat H$ an arbitrary Hilbert space. Now assume
$A_1,A_2$ unital. Note that (by considering e.g.\
$A_1\otimes_{\rm min} A_2)$ we know that there exists a
pair $\pi_i\colon \ A_i\to B({\cal H})$ of unital faithful
representations on the same Hilbert space ${\cal H}$.

\n Consider $\sigma_1,\sigma_2$ and $x$ as above. We need
to show that (6) can be rewritten with unital
representations. Let $p = \sigma_1(1)$ and $q =
\sigma_2(1)$. Note that by augmenting $\widehat H$ if
necessary we can assume that $(1-p)(\widehat H)$ and
$(1-q)(\widehat H)$ are of the same Hilbertian dimension,
and that they are both isometric to some direct sum of
copies of ${\cal H}$.

\n This allows us to define (using $\pi_1$ and $\pi_2$)
$*$-homomorphisms $\hat \pi_i\colon \ A_i\to B(\widehat H)$
such that $$\hat\pi_1(1) = 1-p\quad \hbox{and}\quad \hat
\pi_2(1) = 1-q.$$ Then we define for $a_i\in A_i$
$$\eqalign{\hat \sigma_1(a_1) &= p\sigma_1(a_1)p + (1-p)
\hat\pi_1(a_1)(1-p)\cr
\hat \sigma_2(a_2) &= q\sigma_2(a_2)q + (1-q) \hat
\pi_2(a_2)(1-q).}$$ We have (note $\sigma_1(a_1)p =
p\sigma_1(a_1) = \sigma_1(a_1)\ldots$) $$\sum
\sigma_1(a^1_i) \sigma_2(a^2_i) = p \sum \hat \sigma_1
(a^1_i) \hat \sigma_2(a^2_i)q$$ hence
$$\left\|\sum \sigma_1(a^1_i) \sigma_2(a^2_i)\right\| \le \left\|\sum \hat
\sigma_1(a^1_i) \hat \sigma_2(a^2_i)\right\|$$
but now $\hat\sigma_1, \hat \sigma_2$ are unital
representations (=$*$-homomorphisms), hence this yields
$$\|x\|_h \le \|x\|_{A_1*A_2}.$$
Since the converse is clear (using (6)), this shows that
$p\colon \ E_1 \otimes _h E_2 \to A_1*A_2$ is isometric.
The proof that it is completely isometric is the same with
``operator coefficients'' instead of scalars, we leave
this to the reader.\medskip

\n  We now turn to part (ii). Consider $x$ in
$E \otimes B(H)$ with $\|x\|_{\rm min} <1$. By (i), $x$
corresponds (via $p$) to an element $y$ in $A_1\otimes
A_2$ with $\|y\|_{(A_1\otimes_h A_2) \otimes_{\rm min}
B(H)} <1$. By Lemma~8, we have $y = y_1\odot y_2$ with
$$y_i \in A_i\otimes B(H)\quad \hbox{and}\quad \|y_i\|_{\rm
min}<1.$$ We can write $$y_i = \sum_k a^k_i \otimes b^k_i$$
with $a^k_i \in A_i, b^k_i \in B(H)$ and
$$x = \sum_{k,\ell} a^k_1 a^\ell_2 \otimes b^k_1b^\ell_2.$$
Now consider an isometric representation
$$\pi\colon \ (A_1*A_2) \otimes_{\rm max} B(H) \to B({\cal H})$$
and let $\sigma_1,\sigma_2$ and $\rho$ be its restrictions respectively to
$A_1\otimes 1$, $A_2 \otimes 1$ and $1\otimes B(H)$.
We have (since the ranges of $\rho$ and $\sigma_2$ commute)
$$\eqalign{\pi(x) &= \sum_{k,\ell} \sigma_1(a^k_1) \sigma_2(a^\ell_2)
\rho(b^k_1) \rho(b^\ell_2)\cr
&= \left(\sum_k \sigma_1(a^k_1) \rho(b^k_1)\right) \left(\sum_\ell
\sigma_2(a^\ell_2) \rho(b^\ell_2)\right)\cr
&= \pi(y_1) \pi(y_2).}$$
Hence we conclude that
$$\eqalignno{\|x\|_{\rm max} = \|\pi(x)\| &\le \|\pi(y_1)\| \,
\|\pi(y_2)\|\cr
&\le \|y_1\|_{\rm max} \|y_2\|_{\rm max}\cr
\noalign{\hbox{hence by our assumption on $A_1$ and $A_2$}}
&\le \|y_1\|_{\rm min} \|y_2\|_{\rm min}<1.}$$
This shows by homogeneity 
that $\|x\|_{\rm max} \le\|x\|_{\rm min}$. 

\n Finally, to prove (ii) in the non-unital case, we
simply replace $A_1,\ A_2$ by their unitizations, which
clearly
still satisfy (3). Then, by the unital
case,  $A_1\dot{*}A_2$
appears (see the next remark) as an ideal
in a unital $C^*$-algebra  $A$ such that 
$A \otimes_{\rm min} B(H) = A\otimes_{\rm max} B(H).$
But, as is classical, this property
is inherited by closed ideals (since if $I$ is an ideal
in $A$ and $B$ is
any  other $C^*$-algebra, then the inclusion
$I\otimes_{\rm max} B \subset A\otimes_{\rm max} B$ is
isometric).
\qed

\n {\bf Remark.} Let us denote by $\tilde A$ the unitization of a
$C^*$-algebra $A$. Let $(A_i)_{i\in I}$ be a family of $C^*$-algebras. Then
it is easy to check that the unitization of $\mathop{\dot *}\limits_{i\in
I} A_i$ can be identified canonically with $\mathop{*}\limits_{i\in I}
\tilde A_i$, in short we have
$$\widetilde{\mathop{\dot *}\limits_{i\in I}A_i} \simeq
\mathop{*}\limits_{i\in I} \tilde A_i.$$
Nevertheless, we do not see how to deduce from this
the passage from the non-unital case
(treated in [CES]) to the unital one,
in the first part of the preceding lemma.

\n {\bf Proof of Theorem 7.} It clearly suffices to prove (4) in case $I$
is finite, hence by iteration we may as well assume that $I = \{1,2\}$. Let
$E$ be as in Lemma~9. We will apply Theorem~1 to the
subspace $E\otimes B(H) \subset (A_1*A_2) \otimes B(H)$.
By Lemma~9, the assertion (i) in Theorem~1 is satisfied in
this case (with $E_1,E_2$ now replaced by $E,B(H)$).
Hence, by Theorem~1, we have (4).\qed

Let $C,A$ be $C^*$-algebras. We will denote by $CP(C,A)$ the set of all
completely positive (in short c.p.) maps from $C$ to $A$. A linear map
$u\colon \ C\mapsto A$ is called decomposable if it is a linear combination
of completely positive maps, i.e.\ it can be written as $u = u_1 - u_2	+
i(u_3-u_4)$ with all $u_i$'s completely positive. We denote by $D(C,A)$ the
set of all such maps. This set could be normed by defining for instance
$\|u\| = \inf\left\{ \sum\limits^4_{i=1} \|u_j\|\right\}$ but such a
definition
is not consistent with the algebraic context. Instead, we will use the
following definition due to Haagerup:\ we define $\|u\|_{\rm dec}$ as the
smallest $\lambda\ge 0$ such that there exist $S_1,S_2$ in $CP(C,A)$ such
that $\|S_i\|\le \lambda$, $i=1,2$, and such that
$$x \to \left(\matrix{S_1(x)&u(x^*)^*\cr\cr
u(x)&S_2(x)\cr}\right)$$
is a completely positive map from $C$ to $M_2(A)$. If $u$ is not
decomposable, we set $\|u\|_{\rm dec} =~\infty$.

Haagerup [H1] proved that, equipped with this norm,
 ${D}(C,A)$
becomes a Banach  space. Moreover, he proved $\|u\|_{cb} \le \|u\|_{\rm
dec}$ with equality when $u$ is c.p. Also, if $u$ is self-adjoint, then
$$\|u\|_{\rm dec} = \inf\{\|u_1+u_2\|\mid u=u_1-u_2,\quad u_i\in
CP(C,A)\}.$$
The reader should recall (cf.\ [Pa1]) that $\|u\| = \|u\|_{cb}$ (resp.\ $=
\|u(1)\|$) for any c.p.\ map $u\colon \ C\to A$ (resp.\ when $C$ is assumed
unital), and also that, if $C$ is Abelian, then 
a map $u\colon\ C\to A$ is $c.p.$ iff it is
{\it positive} in the usual sense (=positivity preserving).

Curiously, the dec norm admits several slightly different descriptions.
We start with the most convenient one.

\proclaim Lemma 10. Let $x_1,\ldots, x_n$ be elements in a $C^*$-algebra
$A$ and let $u\colon \ \ell^n_\infty\to A$ be the linear map defined by
$u((\alpha_i)) = \sum \alpha_ix_i$. Then
$$\|u\|_{\rm dec} = \inf\left\{\left\|\sum a_ia^*_i\right\|^{1/2}\,
\left\|\sum b^*_ib_i\right\|^{1/2}\right\} \leqno (7)$$
where the infimum runs over all the decompositions
$x_i=a_ib_i$ with $a_i \in A$ and $b_i\in A$ $(i=1,\ldots,
n)$.

\n {\bf Proof.} If $u$ is positive, i.e.\ if $x_i\ge 0$ for all $i$, this
is very easy:\ the optimal decomposition is simply $x_i =  x^{1/2}_i
x^{1/2}_i$.

\n Let us denote temporarily by $|||(x_i)|||$ the right
side of (7). Assume first that $\|u\|_{\rm dec}<~1$. Then
going back to the definition
of the dec-norm, we can find $y_i,z_i\ge 0$ in $A$ with
$y_i,z_i\ge 0$, $\left\|\sum y_i\right\| < 1$, $\left\|\sum z_i\right\|
<1$ and such that
$$t_i = \left(\matrix{y_i&x^*_i\cr\cr x_i&z_i\cr}\right) \ge 0 \quad
\hbox{for all}\quad i.$$
Then we have (rectangular matrix product)
$$x_i = (0,1)\, t_i\, {1\choose 0}.$$
Let $\gamma_i = (0,1) t^{1/2}_i \in M_{1,2}(A)$ and
$\delta_i = t^{1/2}_i {1\choose 0} \in M_{2,1}(A)$. We
have $x_i = \gamma_i\delta_i$ and $$\left\|\sum
\gamma_i\gamma^*_i\right\| = \left\|\sum z_i\right\| <
1, \quad \left\|\sum \delta^*_i\delta_i\right\| =
\left\|\sum y_i\right\| <1.$$ 
Let $\gamma_i  = (c_i,d_i)$ and $\delta_i = {r_i\choose
s_i}$ so that $x_i = c_ir_i + d_is_i$ with
$$\left\|\sum c_ic^*_i + d_id^*_i\right\| < 1 \quad \hbox{and}\quad
\left\|\sum r^*_ir_i + s^*_is_i\right\| < 1.$$
Assume $A$ unital for simplicity. Let $\varepsilon>0$ and let
$$a_i = (c_ic^*_i + d_id^*_i  + \varepsilon 1)^{1/2} \quad \hbox{and}\quad
b_i = a^{-1}_i x_i.$$
We can choose $\varepsilon>0$ small enough so that
$$\left\|\sum a_ia^*_i\right\| < 1.$$
We have then
$b_i = (a^{-1}_ic_i, a^{-1}_id_i) \delta_i$,\  
hence
$b^*_ib_i = \delta^*_i \omega^*_i\omega_i\delta_i$,\   
where
$\omega_i = (a^{-1}_ic_i, a^{-1}_id_i).$
Note that
$$\eqalignno{\omega_i\omega^*_i &= a^{-1}_i (c_ic^*_i +
d_id^*_i)a^{-1}_i\cr
&\le a^{-1}_i (a^2_i)a^{-1}_i = 1\cr
\noalign{\hbox{hence}}
b^*_ib_i &\le \delta^*_i\delta_i = r^*_ir_i + s^*_is_i,}$$
hence $\left\|\sum b^*_ib_i\right\| < 1$, so that we conclude
$|||(x_i)|||<1$.
 By homogeneity, this completes the proof.\qed

As a consequence, we easily derive the following (known) fact.

\proclaim Lemma 11. Let $n\ge 1$. Let $A$ be a $C^*$-algebra. Then we
have an isometric identity
$$D(\ell^n_\infty, A)^{**} = D(\ell^n_\infty, A^{**}).\leqno (8)$$

\n {\bf Remark.} The reader should be warned that the analogous identity
$cb(\ell^n_\infty, A)^{**} = cb(\ell^n_\infty, A^{**})$ fails to be
isometric in general when $n>2$.

\n {\bf Remark 12.} Concerning the max-norm, we will use several times
the following two known basic facts (see
[L, EL]).\medskip
\item{(i)} For any $C^*$-algebras $A,C$, we have a natural isometric
embedding \break $C\otimes_{\rm max} A \to C\otimes_{\rm
max}A^{**}$ (cf.\ e.g. [W, p. 13]).
 \item{(ii)} Let $C,B,A$ be
$C^*$-algebras and let $\varphi\colon \ A\to B$ be a
completely positive contraction. Then
$I_C\otimes\varphi\colon \ C\otimes_{\rm max}A \to
C\otimes_{\rm max}B$ is a completely positive
contraction (cf.\ e.g. [W, p. 11]). \medskip

\n {\bf Proof of Lemma 11.} This statement reduces to the following fact:\
given an $n$-tuple $(x_i)$ in $A^{**}$ we have $|||(x_i)|||\le 1$ iff there
is a net $(x^\alpha_i)$ of $n$-tuples in $A$ with $|||(x^\alpha_i)||| \le
1$ such that $x^\alpha_i \to x_i\,\ \  \sigma(A^{**},A^*)$
for each $i$. The if part is easy and left to the reader.
To prove the only if part assume without loss of
generality that $|||(x_i)|||<1$ so that $x_i = a_ib_i$
with $a_i,b_i$ in $A^{**}$ such that $\left\|\sum
a_ia^*_i\right\|<1$ $\left\|\sum b^*_ib_i\right\|<1$. Let
$a\in M_n(A^{**})$ (resp.\ $b\in M_n (A^{**}))$ be the
matrix admitting $(a_i)$ (resp.\ $(b_i)$) as its first row
(resp.\ column) and zero elsewhere. Let $a^\alpha$ (resp.\
$b^\alpha$) be a net in the unit ball of $M_n(A)$ tending
$\sigma(M_n(A)^{**}, M_n(A)^*)$ to be $a$ (resp.\ $b$). By
Kaplansky's theorem we can even find a net for which the
convergence is in the strong sense. We define $x^\alpha_i
= a^\alpha(1,i) b^\alpha(i,1)$. Clearly
$|||(x^\alpha_i)||| \le 1$ and $x^\alpha_i \to x_i$ in the
$\sigma(A^{**}, A^*)$-sense. (The strong convergence of
$a^\alpha, b^\alpha$ implies the weak convergence of
$x^\alpha$; moreover on bounded sets weak and
$\sigma$-weak convergences coincide.)\qed

We include here the following simple observation, which is
implicit in [H1, Lemma~3.5].

\proclaim Lemma 13. Let $F$ be a free group
 and let $(U_i)_{i\in I}$ be the  family
of free unitary operators in $C^*(F)$ associated to the
generators, to which we add the unit element. Let
$(a_i)_{i\in I}$ be a finitely supported family in a
$C^*$-algebra $A$ and let $a\colon \ \ell_\infty(I)\to A$
be the mapping defined by $a((\alpha_i)_{i\in I}) =
\sum\limits_{i\in I} \alpha_ia_i$. Then we have
$$\left\|\sum_{i\in I}U_i \otimes
a_i\right\|_{C^*(F)\otimes _{\rm max} A} = \|a\|_{\rm
dec}.\leqno (9)$$

\pf 
Let $A$ be a 
 $\sigma$-finite (=countably decomposable)
von~Neumann algebra. 
By classical facts (cf. e. g. [H4]) we may  assume
$A$ realized as a concrete subalgebra $N\subset B(H)$ and
admitting a (cyclic and) separating vector. Then by [H1,
Lemma~3.5] we have $$\|a\|_{\rm dec} =
\sup\left\{\left\|\sum v_ia_i\right\|\right\}\leqno(10)$$ where the
supremum runs over all unitaries $v_1,\ldots,  v_n$ in
$N'$. Equivalently, if we introduce the representation
$$\pi\colon \ C^*(F) \otimes N\longrightarrow N'\otimes N
\subset B(H\otimes H)$$ defined by $\pi(U_i\otimes a) =
v_ia \qquad \forall a\in N$ then we have $$\|a\|_{\rm dec}
= \left\|\pi \left(\sum U_i\otimes a_i\right)\right\|
\leqno (11)$$ hence this immediately implies
$$\|a\|_{\rm dec} \le \left\|\sum U_i\otimes a_i\right\|_{\rm max}. \leqno
(12)$$
This proves (12) in the 
$\sigma$-finite  von~Neumann case. By a standard direct sum
argument, it is easy to extend
(12) to the case of general von Neumann algebras. But, since the
inclusion $C\otimes_{\rm max} A \subset C\otimes_{\rm max}
A^{**}$ is isometric (see Remark~12) and since (8) holds
we obtain (12) for a general $C^*$-algebra as a consequence
of Lemma~11. The converse inequality $$\left\|\sum U_i
\otimes a_i\right\|_{\rm max} \le \|a\|_{\rm dec}, \leqno
(13)$$ is actually proved in [H1, Lemma~3.5]. Indeed, the
same argument used there shows (without any restriction on $A$)
that for any
representation $\rho\colon \ A\to B({\cal H})$ and for any $v_i$ in
$\rho(A)'$ with $\|v_i\|\le 1$ we have $$\left\|\sum
v_i\rho(a_i)\right\|\le \|a\|_{\rm dec}.$$ This implies in
particular $(13)$. This completes the proof of Lemma~13.

\n An alternate proof of $(13)$ can also be deduced from
(7), as in the above proof of Theorem~0.1.\qed

The following  fact
is due to Kirchberg in [K2]
(and I believe  the simple
proof which follows was known to Kirchberg.)

\proclaim Lemma 14. Let $A$ be a $C^*$-algebra. Assume that
$$C^*(F_\infty) \otimes_{\rm min} A = C^*(F_\infty) \otimes_{\rm max} A.
\leqno (14)$$
Then $A$ is WEP.

\n {\bf Proof.} By Lance's results [L], we know that $A$ has WEP iff for
any embedding $A\to B$ of $A$ into a $C^*$-algebra $B$ and for any
$C^*$-algebra $C$ we have an injective (= isometric) morphism
$$C \otimes_{\rm max} A \to C\otimes_{\rm max} B.$$
Note that this is obvious for the min-norms, therefore if $C\otimes_{\rm
max} A = C \otimes _{\rm min} A$, this is certainly true. Hence, by
assumption, this holds whenever $C=C^*(F)$ with $F$ a free (discrete)
group. 

\n Now consider $C$ arbitrary. Let $F$ be a free group
large enough so that there is a surjective
representation $q\colon\ C^*(F)\to C$. Let $J$ be the
kernel of $q$. Then by the exactness properties of the
max-tensor product (see e.g.\ [W]), we have exact
sequences $$\matrix{0&\to&J\otimes_{\rm max}
A&\to&C^*(F)\otimes_{\rm max}A&\to &C\otimes_{\rm
max}A&\to &0\hfill\cr\cr 0&\to&J\otimes_{\rm max}
B&\to&C^*(F)\otimes_{\rm max}B&\to&C \otimes_{\rm
max}B&\to &0\,\, .\hfill\cr}$$ By the first part of the
proof we have an injective (= isometric) inclusion
$$\varphi\colon \ C^*(F) \otimes_{\rm max} A \to
C^*(F)\otimes_{\rm max} B.$$ Moreover, using an
approximate unit in $J$, it is rather easy to check that
if we view $J\otimes_{\rm max} A$ (resp.\ $J\otimes_{\rm
max} B$) as included in $C^*(F) \otimes_{\rm max}A$
(resp.\ $C^*(F) \otimes_{\rm max}B$) then we have
$$\varphi^{-1}(J\otimes_{\rm max} B) = J\otimes_{\rm max}
A. \leqno (15)$$ Now using (15) and chasing diagrams
it is easy to see that the injectivity of $\varphi$
implies that of the natural map $C \otimes_{\rm max} A \to
C\otimes_{\rm max}B$. Thus, by Lance's criterion (as
mentioned above) $A$ is WEP.\qed

\proclaim Lemma 15. Let $C$ be a $C^*$-algebra. If
$$C\otimes_{\min} B(H)=C\otimes_{\max} B(H)$$
then for any WEP $C^*$-algebra $A$
$$C\otimes_{\min} A=C\otimes_{\max} A.\leqno(16)$$

\pf If $A$ is WEP, by definition 
 we have a factorization $A {\buildrel
\varphi
\over \longrightarrow} B(H) {\buildrel \psi\over\longrightarrow} A^{**}$,
of the canonical inclusion map, with completely
positive contractions $\varphi, \ \psi$.  By Theorem 0.1,
we have a contraction $$I_C\otimes \varphi\colon \
C\otimes_{\rm min} A\to C\otimes_{\rm min} B(H) = C
\otimes_{\rm max} B(H).$$ We  follow this  by $I_C \otimes
\psi$ which is contractive from $C \otimes _{\rm max}
B(H)$ to $C\otimes_{\rm max} A^{**}$ by Remark 12 (ii).
Thus we have a contractive inclusion $C\otimes_{\rm min}A
\to C \otimes_{\rm max} A^{**}$ which (using Remark
12 (i)) implies (16). This shows that
if (16) holds with $A=B(H)$, then it holds whenever $A$ is
WEP. \qed

\n {\bf Remark 16.} In his paper [Ki2], Kirchberg shows
that $$C^*(F_\infty)\otimes_{\rm min} A = C^*(F_\infty)
\otimes_{\rm max} A$$ iff $A$ is WEP. It also follows
from Theorem O.1 and the last two lemmas.

In [K2], Kirchberg also proves a general result on the tensor products
$C\otimes N$ when $N$ is an {\it arbitrary\/} von~Neumann algebra. In that
case, (but with $C$ an arbitrary $C^*$-algebra) we can define a $C^*$-norm
$\|~~\|_{\rm nor}$ on $C\otimes N$ as follows
$$\left\|\sum a_i \otimes b_i\right\|_{\rm nor} = \sup\left\{\left\| \sum
\sigma(a_i) \pi(b_i)\right\|\right\}$$
where the supremum runs over all pairs of representations $\sigma\colon \
C\to B({\cal H})$ $\pi\colon \ N\to B({\cal H})$ with commuting ranges and
with $\pi$ {\it normal\/}. We denote by $C\otimes_{\rm nor} N$ the
completion of $C\otimes N$ for this norm. (See [EL] for more information.)

Our method also allows to prove  Kirchberg's theorem on this tensor norm.

\proclaim Theorem 17. Let $F$ be any free group and let  $C = C^*(F)$. Let
$N$ be any von~Neumann algebra. Then
$$C\otimes_{\rm nor} N = C\otimes_{\rm max} N.$$

\n {\bf Proof.} Replacing $N$ by $M_n(N)$ and using Theorem~1, it clearly
suffices to prove that the norms $\|~~\|_{\rm nor}$ and $\|~~\|_{\rm max}$
are equal on $E\otimes N$ when $E$ is the linear span
of the unit and the free unitary generators of $C$.
Then, we argue as in Lemma 13
(first assuming $N$ $\sigma$-finite, then
passing to the general case): so that 
 by [H1, Lemma~3.5] we find, by (10) and
(11), that  for any $t$ in $E\otimes N$ with associated linear map
$T\colon \ E^*\to N$, we have  $$\|T\|_{\rm dec} \le \|t\|_{\rm
nor}$$ hence (see Lemma~13) $\|t\|_{\rm max} \le
\|t\|_{\rm nor}$.\qed

We conclude this paper with an application to the notion of exactness for
$C^*$-algebras. Recall that a $C^*$-algebra
(or  more generally an operator space) $A$ is called
exact (see [K3]) if for any (closed 2-sided) ideal
$I\subset B$ in $C^*$-algebra $B$, we have an isomorphism
$$B/I \otimes_{\rm min} A \approx {B\otimes_{\rm min} A\over I \otimes_{\rm
min} A}.$$
In [Pi], some of Kirchberg's results on exactness are transferred to the
operator space setting. Let $E$ be a finite dimensional operator space, and
let
$$u_E\colon \ B/I \otimes_{\rm min}E \to {B\otimes_{\rm min} E\over I
\otimes_{\rm min} E}$$
be the canonical isomorphism.

Let $F$ be another finite dimensional operator space. Recall the notation
$$d_{cb}(E,F) = \inf\{\|u\|_{cb} \|u^{-1}\|_{cb}\}$$
where the infimum runs over all possible isomorphisms $u\colon \ E\to F$.
By convention we set $d_{cb}(E,F) = \infty$ if $E,F$ are not isomorphic.

We denote
$$d_{SK}(E) = \inf\{d_{cb}(E,F)\mid F\subset K(\ell_2)\}.$$
(Here $K(\ell_2)$ denotes the algebra of all compact operators on
$\ell_2$.)

When $\widehat E$ is an infinite dimensional operator space, we define
$$d_{SK}(\widehat E) = \sup\{d_{SK}(E) \mid E\subset \widehat E \quad
\dim(E) < \infty\}.$$
In [Pi], by a simple adaptation of an argument of Kirchberg in [K3], we
show that for any exact operator space $E$
$$d_{SK}(E) = \sup\{\|u_E\|\}= \sup\{\|u_E\|_{cb}\}\leqno (17)$$
where the supremum runs over all possible pairs $(I,B)$ with $I\subset B$.
(Actually, it suffices to consider $I = K(\ell_2)$ and $B = B(\ell_2)$). In
[K3], Kirchberg showed that a $C^*$-algebra $A$ is exact iff $d_{SK}(A)
=1$. The point of the next result is that it suffices for the exactness of
$A$ to be able to embed (almost completely isometrically) the linear span
of the unitary generators of $A$ and the unit into $K(\ell_2)$ (or into a
nuclear $C^*$-algebra).

\proclaim Theorem 18. Let $E\subset A$ be a closed subspace of a unital
$C^*$-algebra $A$. We assume that $1_A\in E$ and that $E$ is the closed
linear span of a family of unitary elements of $A$. Moreover, we assume that
$E$ generates $A$ (i.e.\ that the smallest $C^*$-subalgebra of $A$
containing $E$ is $A$ itself). Then, if $d_{SK}(E)=1$, $A$ is exact.

\pf Let $(I,B)$ be as above with $B$ unital. By (17), if $d_{SK}(E)=1$, the
unital $*$-homomorphism
$$\pi\colon \ B/I\otimes A \longrightarrow {B\otimes_{\rm min} A\over
I\otimes_{\rm min}A}$$
becomes completely contractive when restricted to $(B/I) \otimes_{\rm min}
E$. By Proposition~6, $\pi$ extends to a continuous (= contractive)
$*$-homomorphism on $(B/I) \otimes_{\rm min}A$. Hence $A$ is exact. \qed
\vfill\eject

\centerline{\bf References}

\item{[BP]} D. Blecher and V. Paulsen. Tensor products of operator spaces.
 J. Funct. Anal. 99 (1991) 262-292.
 
\item{[Ch1]} M. D. Choi. A Schwarz inequality for positive linear maps on
$C^*$-algebras. Illinois J. Math. 18 (1974) 565-574.

\item{[Ch2]} M. D. Choi.  A simple $C^*$-algebra
 generated by two
finite order unitaries. Can. J. Math. 31 (1979) 887-890.

\item{[CS]}  E.  Christensen and A. Sinclair.  A survey
of completely bounded operators.	Bull. London Math. Soc.
21 (1989) 417-448.

\item{[CES]}  E.  Christensen, E. Effros and A. Sinclair. 
Completely bounded multilinear maps and
$C^*$-algebraic cohomology.
Invent. Math. 90 (1987) 279-296.

 \item{[EL]} E. Effros and C. Lance. Tensor products of
operator algebras.  Adv. Math. 25 (1977) 1-33.

 \item{[ER1]} E. Effros and Z.J. Ruan. On matricially normed spaces. Pacific
 J. Math. 132 (1988) 243-264.

 \item{[ER2]} E. Effros and Z.J. Ruan.  Self duality for
the Haagerup
 tensor product and Hilbert space factorization.  J. Funct. Anal. 100
(1991) 257-284.

\item{[H1]} U. Haagerup. Injectivity and decomposition of completely
 bounded maps in ``Operator algebras and their connection with Topology and
 Ergodic Theory''. Springer Lecture Notes in Math. 1132 (1985) 91-116.
 
 \item{[H2]} $\underline{\hskip1.5in}$. An example of a non-nuclear
 $C^*$-algebra which
 has the metric approximation property. Inventiones Math. 50 (1979) 279-293.
 
 \item{[H3]} $\underline{\hskip1.5in}$. Decomposition of completely bounded
 maps on operator algebras. Unpublished manuscript. Sept. 1980.
 
 \item{[H4]} U. Haagerup. The standard form of a von
Neumann algebra.
Math. Scand. 37 (1975) 

\item{[Ha]} A. Harcharras. On some stability properties
of the full $C^*$-algebra associated to the free group
$F_\infty$. In preparation.

 \item{[K1]} E. Kirchberg. On non-semisplit
extensions, tensor products and exactness of group
$C^*$-algebras. Invent. Math. 112 (1993) 449-489.

 \item{[K2]} E. Kirchberg. Commutants of
unitaries in UHF algebras and functorial properties of
exactness. J.
 reine angew. Math. 452 (1994)
39-77. 
 
\item{[K3]} E. Kirchberg. On subalgebras of the
CAR-algebra. 
J. Funct. Anal. 129 (1995) 35-63.

 \item{[L]} C. Lance. On nuclear $C^*$-algebras. J. Funct.
Anal. 12 (1973) 157-176.

 \item{[Pa1]} V. Paulsen.  Completely
bounded maps and dilations. Pitman Research Notes 146.
Pitman Longman (Wiley) 1986.
 
\item{[Pa2]}  V. Paulsen. Representation of
Function algebras, Abstract
 operator spaces and Banach space Geometry. J.
Funct. Anal. 109 (1992) 113-129.
 
\item{[Pa3]}  V. Paulsen. The maximal
operator space of a normed space. Preprint to appear.
 
 \item{[PS]} V. Paulsen and R.Smith. Multilinear maps and
tensor norms on operator systems. J. Funct. Anal. 73
(1987) 258-276.

\item{[Pi]}  G. Pisier. Exact operator spaces. Proceedings Colloque sur
les alg{\`e}bres d'op{\'e}rateurs (Orl{\'e}ans 1992)
 Ast{\'e}risque, Soc. Math. France (To
appear in 1995).

\item{[R]} M. Rieffel. Induced representations of
$C^*$-algebras. Adv. Math. 13 (1974) 176-257.

\item{[W]} S. Wassermann. Exact $C^*$-algebras and
related topics.
Lecture Notes Series, Seoul Nat. Univ. (1994).
 
\vskip12pt

Texas A\&M University

College Station, TX 77843, U. S. A.

and

Universit\'e Paris VI

Equipe d'Analyse, Case 186,
 
75252 Paris Cedex 05, France

 \end